\documentclass[twoside,12pt]{amsart}

\usepackage{amsmath}
\usepackage{amsthm}

\numberwithin{equation}{section}

\input epsf 

\def\deg {\hbox{deg}}
\def\ll#1{\lambda_#1}
\def\pder#1#2{\frac {\partial#2} {\partial#1}}
\def\pderr#1#2#3{\frac {\partial^#2#3} {\partial#1^#2}}
\def\spder#1#2#3{\frac {\partial^2#3} {\partial#1\, \partial#2}}
\def\ul{\frac {u^{\lambda}} {\lambda!}}
\def\deg {\hbox{deg}}
\def\Rbar{\bar R}

\def\figureform#1#2{
\bigskip
\centerline{\epsfbox{#1}}

\centerline{#2}
\bigskip}

\def\set#1{\{\,#1\,\}}

\def\Tv{\widetilde T}
\def\Rv{\widetilde R}

\newtheorem*{lem}{Lemma}

\def\smbx{\vbox{\hrule\hbox{\vrule\kern3pt
		\vbox{\kern6pt}\kern3pt\vrule}\hrule}}

\begin{document}

\title[Hypergraphs and a functional equation]{Hypergraphs and a functional equation of Bouwkamp and de Bruijn}

\author{Ira M.~Gessel}
\address{ Department of Mathematics\\ Brandeis University\\ Waltham, MA 02254-9110}
\email{gessel@brandeis.edu}
\author{Louis H.~Kalikow}
\address{Department of Mathematics\\ Haverford College \\Haverford, PA 19041-1392}
\email{lkalikow@haverford.edu}

\begin{abstract}
Let $\Phi(u,v)=\sum_{m=0}^\infty  \sum_{n=0}^\infty c_{mn}u^m v^n$. Bouwkamp and de~Bruijn found that there exists a
 power series $\Psi(u,v)$ satisfying the equation $t\Psi(tz,z)=\log \biggl(\sum_{k=0}^\infty  \frac {t^k} {k!} \exp(k\Phi(kz,z))\biggr) $. We show that this result  can be interpreted combinatorially using hypergraphs. We also explain some facts about $\Phi(u,0)$ and $\Psi(u,0)$, shown by Bouwkamp and de~Bruijn, by using hypertrees, and we use Lagrange inversion to count hypertrees by number of vertices and number of edges of a specified size.
\end{abstract}
\maketitle

\section{Introduction}
\label{origsec1}
In \cite{BB}, Bouwkamp and de~Bruijn use algebraic
methods to prove some results concerning a power series 
expansion. Their original motivation arose from work by Harris and Park \cite{HP}, who showed the asymptotic normality of the distribution of empty cells when some number of balls were placed in some number of equiprobable cells. To accomplish this, Harris and Park employed factorial cumulants; in particular, they showed that in
\begin{equation*}
\log\biggl(\sum_{k=0}^{\infty} \binom N k  \left(1-\frac k N\right)^N t^k\biggr),
\end{equation*}
the coefficient of $t^n$ is $O(N)$. As noted in \cite{BB} and \cite{HP}, de Bruijn did some work on this problem, and it led  Bouwkamp and de~Bruijn to show that if $\Phi(u,v)$ is a double 
power series of the form
\begin{equation*}
\Phi(u,v)=\sum_{m=0}^\infty  \sum_{n=0}^\infty c_{mn}u^m v^n,
\end{equation*}
then there exists a power series $\Psi(u,v)$ such that
\begin{equation}
\log \biggl(\sum_{k=0}^\infty  \frac {t^k} {k!} \exp(k\Phi(kz,z))\biggr)= 
t\Psi(tz,z).\label{orig1}
\end{equation}
That is, the left side can be written as $t\sum_{n=0}^{\infty} t^n 
\theta_n (z)=\sum_{n=0}^{\infty} t^{n+1} 
\theta_n (z)$, 
where $\theta_n(z)$ is a power series which has no powers of $z$ less than $n$. 

Bouwkamp and de~Bruijn further 
demonstrate a result relating $\psi(u):=\Psi(u,0)$ and $\phi(u):=\Phi(u,0)$. Note 
that $\Psi(u,0)$ 
yields the  ``leading 
terms'' of \eqref{orig1}, in the sense that $\,t\Psi(tz,0)$ is the series which contains all 
terms of 
$\,t\Psi(tz,z)$ in which the power of $z$ is one less then the power of $t$.  
Bouwkamp and de~Bruijn show 
that if $w$ is the power series in $y$ satisfying
\begin{equation}
y=w\exp (-\phi(w)-w\phi'(w)),\label{orig2}
\end{equation}
then 
\begin{equation}
\psi(y)=(w-w^2\phi'(w))/y.\label{orig3}
\end{equation}

We will show that these results from \cite{BB} are actually consequences of identities 
for hypergraphs and hypertrees. We will also give combinatorial interpretations of many other equations that were derived algebraically in [3].

A hypergraph is a generalization of a graph (the next section has exact definitions and basic facts;
see \cite{B} for further background). 
In general, edges 
can consist not only of a set of two vertices, but of a set of an arbitrary number 
of vertices. An edge consisting of $i$ vertices will be called an 
{\it $i$-edge\/}. 
We will be concerned with hypergraphs without empty edges or loops 
(i.e., without 0-edges or 1-edges); therefore, when we use the term {\it hypergraph\/}, it will refer 
to hypergraphs whose edges have at least two vertices.

Bouwkamp and de~Bruijn prove their results by analyzing the power
series (see (1.2) in \cite{BB})
\begin{equation}
\sum_{k=0}^{\infty} \frac {t^k} {k!} \exp \biggl[\sum_{i=2}^{\infty} k^i
x_{i-1}\biggr]\label{orig5}
\end{equation}
and then substituting power series for $t$ and the $x_i$.
We shall prove their results in a very similar way, by substituting power 
series for $t$ and for the $u_i$ in the generating function
\begin{equation}
\sum_{k=0}^{\infty} \frac {t^k} {k!} \exp \biggl[\sum_{i=2}^{\infty} \binom k i u_i\biggr],\label{orig4}
\end{equation}
which we interpret as a generating function for hypergraphs.

If we wanted, we could prove the result by considering \eqref{orig5} to be the 
exponential 
generating 
function for a set of  objects in which the ``edges'' are sequences of vertices (with repetitions allowed). However, since hypergraphs seem more natural, we 
choose to use 
them.

Section 2 provides definitions of hypergraphs and hypertrees. In section \ref{origsec3}, we 
prove \eqref{orig1} by showing it is a consequence of the hypergraph analogue of the fact that every connected graph with $n$ vertices has at least $n-1$ edges.
 In Section 4, we interpret several equations obtained by
Bouwkamp and de~Bruijn in terms of hypergraphs. Section \ref{origsec4} provides 
interpretations of the 
leading terms, $\psi(u)$ and $\phi(u)$, using hypertrees. We conclude in Section 6 
by showing how this work and Lagrange inversion can be used to obtain previously known results on the enumeration of hypertrees.

\section{Definitions and background}
\label{origsec2}
We define a hypergraph $H$ on $n$ vertices to be an ordered pair $(V,E)$, where 
$V$ is the set of vertices, with $|V|=n$, and $E$ is a multiset of subsets
of $V$; we also require the the subsets in $E$ contain at least two vertices. In particular, 
we allow 
multiple edges. For an arbitrary hypergraph $H$, we let $v(H)$ denote the number 
of vertices of 
$H$ and $e(H)$ denote the number of edges of $H$.  This definition differs from that in Berge \cite{B} 
since we allow a hypergraph to have vertices which belong to no edge.
Our definition of a hypergraph nearly agrees
with that of 
Grieser \cite{Gr}; the difference is that we do not allow loops.

In general, we will consider hypergraphs labeled so that if the hypergraph has 
$n$ 
vertices, they are labeled by the elements of $[n]:=\set{1,2,3,\ldots,n}$, and if the hypergraph has 
$\lambda_i$ 
$i$-edges, they are labeled by the elements of $[\lambda_i]$.
For simplicity, we will call such objects {\it labeled hypergraphs\/}.

In what follows, we will always have $\ll 1=0$, since our hypergraphs have
no loops. Let $u_2, u_3, u_4, \ldots$ be indeterminates. We  define the 
{\it weight\/} of $H$ to be
\begin{equation*}
u_2^{\lambda_2}u_3^{\lambda_3}  \cdots u_{n}^{\lambda_{n}},
\end{equation*}
and we define the {\it edge magnitude\/} of $H$ to be $\sum_{i=2}^{k} (i-
1)\lambda_i$. 

An example of a labeled hypergraph is given in Figure~1. The 2-edges are denoted by a 
segment 
connecting the two vertices; for edges with 
more than two 
vertices, the edge is 
represented by a closed curve which contains the vertices of the edge inside it. 
The vertices are 
labeled by numbers without
 subscripts; for clarity, the edges are labeled with subscripted numbers in which the subscript refers to the size of the edge being labeled. (The subscripts on the edge labels thus do not add structure to the hypergraph.)
 For the hypergraph in the figure, $V=[5]$; 
$E=\{\,\set{1,2}$,
$\set{1,2}$, $\set{3,5}$, $\set{4,5}$, $\set{1,3,4}$, $\set{1,3,4}$, 
$\set{3,4,5}\,\}$; $\ll 2=4, \ll 3=3$; the weight  is 
$u_2^4 u_3^3$; and the edge magnitude is $10$.

\figureform{bbfig.1}{Figure 1: A sample labeled hypergraph}
\bigskip

We define a {\it walk\/} in a hypergraph to be a sequence
\begin{equation*}
v_0,e_1,v_1, \ldots, v_{n-1}, e_n, v_n,
\end{equation*}
where for all $i$, $v_i \in V$, $e_i \in E$, and
for 
each $e_i$, $\{v_{i-1},v_i\} 
\subseteq e_i$.
We define a {\it path\/} in a hypergraph to be a 
walk in which 
all $v_i$ are distinct  and all $e_i$ are distinct. A walk is a {\it cycle\/} if the walk 
contains at 
least two edges, all $e_i$ are distinct, and all $v_i$ are 
distinct except $v_0=v_n$.

A hypergraph is {\it connected\/} if for every pair of vertices $v,v'$ in the 
hypergraph, there is a path starting at $v$ and ending at $v'$. The hypergraph  in 
Figure~1 is 
connected. For example, a path between vertices 2 and 5 is

\begin{equation*}
2,\, 1_2=\set{1,2},\, 1,\, 3_3=\set{1,3,4},\, 3,\, 3_2=\set{3,5},\, 5.
\end{equation*}
We define a {\it hypertree\/} to be a connected hypergraph with no cycles.  

The {\it degree\/} of a vertex $v \in V$, denoted $\deg (v)$, is defined as 
\begin{equation*}
\deg (v):= |\{e\in E \mid v \in e\}|;
\end{equation*}
i.e., the degree of $v$ is the number of edges to which $v$ belongs. Two 
vertices in a 
hypergraph are {\it adjacent\/} if there is an edge containing both. 

We now note some basic facts about hypertrees. First, two edges in a hypertree 
have at most one vertex in common; for if 
edges $e_1, e_2$ have two vertices $v_1, v_2$ in common, then the hypergraph 
has a cycle $v_1,e_1,v_2,e_2,v_1$. Next, we prove the following lemma
It is known; for another proof, see \cite{Gr}.

\begin{lem} 
A connected hypergraph on $n$ vertices is a hypertree if and only if it has edge magnitude $n-1$. Furthermore, the minimum edge magnitude of a connected hypergraph is $n-1$.
\end{lem}

\begin{proof}

First we prove by induction on $n$ that a hypertree on $n$ vertices
has edge magnitude $n-1$. This is clearly true for $n=1$. Now suppose
that $H$ is a hypertree with $n>1$ vertices and that every hypertree
with $n-1$ vertices has edge magnitude $n-2$.

Let $v_0,e_1,\dots, e_m,v_m$ be a longest path in $H$. Suppose that
$v_m$ is contained in some edge $e$ other than $e_m$. We will show
that this assumption leads to a contradiction. If $e\in\{e_1,\dots,
e_{m-1}\}$ then $H$ contains a cycle, so $e\notin\{e_1,\dots,
e_{m-1}\}$. Let $v$ be vertex of $e$ other than $v_m$. Then
$v_0,e_1,\dots, e_m,v_m,e,v$ is either a longer path than
$v_0,e_1,\dots, e_m,v_m$ (if $v\notin\{v_0,\dots,v_{m-1}\}$) or
contains a cycle (if $v\in\{v_0,\dots,v_{m-1}\}$). Since both are
impossible, $v_m$ cannot be contained in any edge of $H$ other than
$e_m$.

Let $H'$ be obtained from $H$ by removing vertex $v_m$ and then either replacing
edge $e_m$ with 
$e_m - \{v_m\}$, if $|e_m|>2$; or deleting $e_m$, if $|e_m|=2$. It is
clear that $H'$ is a
hypertree and that its edge magnitude is one less than that of $H$. By
the inductive
hypothesis, $H'$ has edge magnitude $n-2$, so $H$ has edge magnitude
$n-1$.

Now suppose that $H$ is a connected hypergraph on $n$ vertices that is
not a hypertree. Then $H$ has a cycle $v_0, e_1, v_1,\dots, e_n,
v_0$. Replacing $e_1$ with $e_1-\{v_0\}$ if $|e_1|>2$, or deleting
$e_1$ if $|e_1|=2$, leaves a connected hypergraph on $n$ vertices with
edge magnitude one less than that of $H$. Repeating this reduction
eventually yields a hypertree. Thus $H$ has edge magnitude greater
than $n-1$.
\end{proof}

\section {Proof of {\eqref{orig1}}}
\label{origsec3}

We turn our attention to proving the main result of \cite{BB} using the 
exponential generating function for labeled hypergraphs. We adopt 
the convention that if $\lambda=(\lambda_2, \lambda_3, \ldots\,)$ is a sequence of 
integers
with finitely many non-zero parts, then
\begin{equation*}
u^{\lambda}= u_2^{\lambda_2}u_3^{\lambda_3}u_4^{\lambda_4}\cdots,
\end{equation*}
and
\begin{equation*}
\frac {u^{\lambda}} {\lambda !}
 = \frac {u^{\lambda}} {\lambda_2 ! \, \lambda_3 ! \, \lambda_4 ! \, \cdots}.
\end{equation*}
We do not use $u_1$ since we will not consider hypergraphs with loops.

Using exponential generating functions, we now count labeled hypergraphs with 
vertices and 
edges labeled as in Figure~1. (For background on the combinatorics of exponential 
generating functions, see \cite[Chapter 3, Section 2]{GJ}, \cite[Chapter 5]{Stan}, or, for an approach using species, \cite[Chapter 1]{BLL}). Consider the 
exponential generating function
\begin{equation*}
(e^{u_i})^{\binom k i}=(1+ u_i + \frac {u_i^2} {2!} + \cdots)^{\binom k i}.
\end{equation*}
We view the term $\frac {u_i^j} {j!}$ in the expansion of $e^{u_i}$ as representing  
$j$ multiple copies of a particular $i$-edge. Since there are $\binom k i$ 
$i$-subsets of vertices in a
hypergraph with $k$ vertices, the previous expression counts labeled hypergraphs on $k$ 
vertices whose edges are all of size $i$. Therefore,
\begin{equation*}
\prod_{i=2}^{\infty} (e^{u_i})^{\binom k i} 
= \exp \biggl[\binom k 2 u_2 + \binom k 3 u_3 + \cdots\biggr]
\end{equation*}
is the exponential generating function for labeled hypergraphs on $k$ vertices, 
where the 
coefficient of $\ul$ 
is the number of labeled hypergraphs on $k$ vertices with $\ll i$ 
$i$-edges for each  $i\ge 2$. 
From this, we see that the exponential generating function for labeled hypergraphs 
with vertices weighted $t$ and $i$-edges weighted by $u_i$ is
\begin{equation*}
\sum_{k=0}^{\infty} \frac {t^k} {k!} \exp \biggl[\sum_{i=2}^{\infty} {\binom k i}u_i\biggr].
\end{equation*}
Since the edge magnitude of a hypergraph counted by the 
coefficient of $\ul$ is $\sum_{i\ge 2} \ll i (i-1)$, we define the {\it magnitude\/} 
of 
$u^{\lambda}$
to be the same expression.

We now consider connected labeled hypergraphs. Since a 
labeled hypergraph 
is a set of connected labeled hypergraphs, if $C:=C(t;u_2,u_3,\ldots)$ is the exponential 
generating 
function for connected labeled hypergraphs, we have 
\begin{equation*}
e^C=\sum_{k=0}^{\infty} \frac {t^k} {k!} 
\exp \biggl[\sum_{i=2}^{\infty} {\binom k i}u_i\biggr].
\end{equation*}
Hence,
\begin{equation}
C=\log\biggl[\sum_{k=0}^{\infty} \frac {t^k} {k!} 
\exp \biggl({\binom k 2} u_2 + {\binom k 3} u_3 + {\binom k 4} u_4 + 
\cdots\biggr)\biggr].\label{orig3.1}
\end{equation}

We know from  Section~2 that the edge magnitude of a connected hypergraph on 
$k$ vertices
 must be at least $k-1$. 
So if we write
\begin{equation}
C=\sum_{k=1}^{\infty} \frac {t^k} {k!} f_k(u_2,u_3,u_4,\ldots),\label{orig3.2}
\end{equation}
the minimum magnitude of terms of $f_k$ is $k-1$.

We finish the proof of \eqref{orig1} with an argument from Bouwkamp and de~Bruijn 
\cite[Section~1]{BB}. For $m\ge 1$, let 
$P_m(z)=\sum_{i\ge 0} p_{m,i} z^i$ be power series in $z$. If we make the substitutions
\begin{equation*}
 t\exp(P_1(z))\mapsto t,\hskip 1in  z^{m-1}P_m(z)\mapsto u_m,
\end{equation*}
in \eqref{orig3.2}, then the coefficient of $t^k$ is a power series in $z$ with no term of 
degree less 
than $k-1$. Thus, the resulting power series is of form $t\Psi(tz,z)$. If we make the same substitutions in \eqref{orig3.1}, we obtain
\begin{equation*}
C=\log\left[\sum_{k=0}^\infty \frac{t^k}{k!} \exp\biggl(
\sum_{m\ge1} {\binom k m} z^{m-1}\sum_{i\ge0} p_{m,i} z^i\biggr)\right].
\end{equation*}
Here the argument of the exponential function may be written
\begin{equation}
\label{step1}\sum_{j=0}^\infty z^j \sum_{m=1}^{j+1} {\binom k m} p_{m,j-m+1}. 
\end{equation}
Now let $\Phi(u,v)=\sum_{m,n=0}^\infty c_{mn}u^mv^n$ be any double power series, so
\begin{equation}\label{step2}
k\Phi(kz, z)=\sum_{m,n\ge0} k^{m+1} c_{mn}z^{m+n}
  =\sum_{j=0}^\infty z^j\sum_{l=1}^{j+1} k^l c_{l-1, j-l+1}.
\end{equation}
Since the polynomials ${\binom k m}$, for $m=1,\dots, j+1$, form a basis for the
polynomials of degree at most $j+1$ in $k$ that vanish at 0, it is possible to
choose  power series $P_m(z)$ that make \eqref{step1} equal to \eqref{step2} and thus make \eqref{orig3.1} equal to the left side of \eqref{orig1}. This proves \eqref{orig1}.

\section {Further combinatorial interpretations}
\label{newfurther}

We remark here that some calculations done by Bouwkamp and de~Bruijn \cite[Sections~2--3]{BB} 
correspond to simple manipulations involving $C$ which yield various ways of 
decomposing hypergraphs into other hypergraphs. Differentiation of $e^C$ with 
respect to $u_j$ yields (cf.~\cite[Section~2]{BB})

\bigbreak
\bigskip
\centerline{\hfil\epsfbox{bbmorefig.5}\hfil\hfil\epsfbox{bbmorefig.6}\hfil}
\smallskip
\centerline{(a) Removing the 2-edge yields two hypergraphs}
\bigskip
\bigskip
\centerline{\hfil\epsfbox{bbmorefig.7}\hfil\hfil\epsfbox{bbmorefig.8}\hfil}
\smallskip
\centerline{(b) Removing the 2-edge yields one hypergraph}
\bigskip
\centerline{Figure 2: Decomposition of \eqref{orig3.4}}
\bigskip\bigskip

\begin{equation}\label{orig3.3}
\begin{split}
\pder {u_j} {e^C} &= \sum_{k=0}^{\infty} \frac {t^k} {k!} {\binom k j}
\exp \biggl(\binom k 2 u_2 + {\binom k 3} u_3 + {\binom k 4} u_4 + 
\cdots\biggr)\\
&=\frac {t^j} {j!} \pderr t j {e^C}.
\end{split}
\end{equation}

By properties of exponential generating functions \cite[pp.~167--8]{GJ}, $\pder {u_j} {e^C}$
counts hypergraphs rooted at an unlabeled $j$-edge. (The generating function $u_j\pder {u_j} {e^C}$ would count those rooted at a labeled $j$-edge.)
Also, operating on $e^C$ by $t^j 
\pderr t j {}$ counts hypergraphs which are equipped with an ordered $j$-tuple of $j$ 
distinct 
vertices. By dividing by $j!$ as in \eqref{orig3.3}, we count hypergraphs rooted at $j$ 
vertices. 
Therefore \eqref{orig3.3}  represents a bijection between hypergraphs 
rooted at an unlabeled $j$-edge and hypergraphs with $j$ rooted vertices.

From  \eqref{orig3.3} in the case $j=2$, we obtain (cf.~\cite[(2.1)]{BB})
\begin{equation}
\pder {u_2} C = \frac {t^2} {2!} \biggl[\biggl(\pder t C\biggr)^2 + \pderr t 2 
C\biggr].\label{orig3.4}
\end{equation}
This represents a way to decompose connected hypergraphs rooted at an unlabeled 2-edge. By 
removing the 
rooted 2-edge, we either obtain two vertex-rooted connected hypergraphs or else 
we obtain a 
doubly-vertex-rooted connected hypergraph. See Figure~2 for an example; in the figure, rooted objects are marked heavily (thick lines or larger dots). Since $(t \pder t C)^2$ 
counts 
ordered pairs of rooted hypergraphs, $\frac 1 {2!} (t \pder t C)^2$ 
counts sets containing two rooted 
hypergraphs, as in Figure 2(a). Also, since $t^2 \pderr t 2 C$ counts hypergraphs rooted at an ordered 
pair of vertices, $ \frac {t^2} {2!} \pderr t 2 C$ counts doubly-vertex-rooted 
hypergraphs, as in Figure 2(b).

\bigskip

Next, note that \eqref{orig3.3} implies 
\begin{equation*}
\pder {u_j} {e^C}=\frac 1 j \biggl(t \pder t {} -(j-1)\biggr) \pder {u_{j-1}} {e^C}.
\end{equation*}

From this we obtain (cf.~\cite[(2.2)]{BB})
\begin{equation}
\pder {u_j} C=\frac 1 j \biggl[t\pder t C \pder {u_{j-1}} C+t \spder t {u_{j-1}} 
C - (j-1) \pder {u_{j-1}} C\biggr].\label{orig3.5}
\end{equation}

\begin{tabular}{cc}
&   \epsfbox{bbmorefig.2}  \\
& (i) \\
 \epsfbox{bbmorefig.1} & \epsfbox{bbmorefig.3}\\
(a) Original hypergraph & (ii)\\
&  \epsfbox{bbmorefig.4}\\
& (iii) \\
& \\
&(b) Decompositions\\
\end{tabular}

\bigskip
\centerline{Figure 3: Decompositions according to \eqref{orig3.5}}
\bigskip\bigskip

Combinatorially, 
the term $t\pder t C \pder {u_{j-1}} C$ on the right side of $\eqref{orig3.5}$ counts 
pairs of connected hypergraphs, in which one of the pair is  rooted at a 
vertex and one is rooted at an unlabeled $(j-1)$-edge. The next term, 
$t \spder t {u_{j-1}} C$, counts connected 
hypergraphs rooted at both an unlabeled $(j-1)$-edge and a vertex. Finally, $(j-1) \pder {u_{j-
1}} C$  counts connected hypergraphs rooted at an unlabeled $(j-1)$-edge 
and at a vertex {\it in that rooted edge\/}. Thus,  \eqref{orig3.5} 
says there are $j$ ways to decompose a connected hypergraph rooted at an unlabeled $j$-edge 
either into a pair of connected hypergraphs, one rooted at a vertex and the other 
rooted at an unlabeled $(j-1)$-edge (for example, see Figure~3(b)(i)); or into a single connected hypergraph, rooted at an unlabeled 
$(j-1)$-edge and at a vertex {\it not in the rooted edge\/} (see Figure~3(b)(ii), (iii)). 
To perform the decomposition given a hypergraph rooted at an unlabeled $j$-edge, simply choose a vertex $v$ in the rooted edge $e$. Then remove the edge $e$, add the edge $e-\set{v}$, and root the 
new object at $v$ and at the added edge. 
\section {Interpretations of the Leading Terms}
\label{origsec4}

We now consider the combinatorial interpretation of the results in \cite{BB} about the 
leading 
terms of \eqref{orig1}. It turns out that much of the work leading to the results in \cite{BB} 
involves 
differential equations related to decompositions of hypertrees.

We define for $n\ge 1$, 
\begin{equation*}
C_n=C_n(u_2,u_3, \ldots):=\biggl[\frac {t^n} {n!}\biggr] C;
\end{equation*}
that is, $C_n$ is the coefficient of $\frac {t^n}{n!}$ in the power series $C$. Thus, the coefficient of $\ul$ in $C_n$ is the number of connected labeled
hypergraphs on 
$[n]$ with $\lambda_j$ $j$-edges for $j\ge 2$. From  \eqref{orig3.4} we obtain 
(cf.~\cite[(2.3)]{BB}),
\begin{equation}
\pder {u_2} {C_n} 
= \frac 1 {2!} \left[\sum_{i=1}^{n-1}{\binom n i} i (n-i) C_i C_{n-
i}+n(n-1)C_n \right],\label{orig4.1}
\end{equation}
and from  \eqref{orig3.5} (cf.~\cite[(2.4)]{BB}),
\begin{equation}
\pder {u_j} {C_n}=\frac 1 {j}\left[
\sum_{i=1}^{n-1} {\binom n i}(n-i) \left(\pder {u_{j-1}} {C_i} \right) C_{n-i} +\bigl(n-(j-1)\bigr)\pder {u_{j-1}} {C_n} \right].\label{orig4.2}
\end{equation}

Now, we define
\begin{equation*}
T_n=T_n(u_2, u_3, \ldots):=\hbox{all terms of magnitude $n-1$ in $C_n$},
\end{equation*}
and let
\begin{equation}\label{orig4.3}
T(t;u_2, u_3, \ldots)=\sum_{n=1}^{\infty} \frac {t^n} {n!} T_n(u_2, u_3, \ldots).
\end{equation}
The generating function $T$ contains the ``leading terms'' of $C$, in the sense $T_n$ is the sum of the terms of minimal 
magnitude in $C_n$.
Hence, if $\tau_{n,\lambda }$ is the number of 
labeled hypertrees on $[n]$ with $\ll i$ $i$-edges, we have
\begin{equation*}
T_n=\sum_{\lambda=(\ll 2, \ll 3, \ldots)}\tau_{n,\lambda }\ul.
\end{equation*}
where $\tau_{n,\lambda} =0$ unless $\sum_{i\ge 2} (i-1)\ll i =n-1$. We note that $T_n$ 
corresponds to $\eta_n^*$ in \cite[Section 3]{BB}. 

We can get differential 
equations for $T$ using the differential 
equations for $C_n$. By taking terms with the minimal magnitude $n-2$ on both sides of 
\eqref{orig4.1} we get (cf.~\cite[(3.1)]{BB})
\begin{equation*}
\pder {u_2} {T_n}= \frac 1 {2!} \biggl(\sum_{i=1}^{n-1} {\binom n i} i(n-i)T_i T_{n-
i}\biggr).
\end{equation*}
There is no contribution  to this equation from the term $\frac 1 2n(n-1)C_n$ in the 
summation on the right side 
of \eqref{orig4.1}. 
That term corresponds to the case in which
the removal of a 2-edge from a connected hypergraph yields a single 
connected hypergraph. For hypertrees, removing a 2-edge must yield two 
hypertrees. 

From the last equation, we obtain
\begin{equation}\label{orig4.4}
\pder {u_2} T=\frac 1 {2!} \biggl(t \pder t T\biggr)^2 .
\end{equation}

If we take terms with the minimal magnitude $n-1-(j-1)=n-j$ on both sides of  \eqref{orig4.2}, we get
(cf.~\cite[(3.2)]{BB})
\begin{equation*}
 \pder {u_j} {T_n}= \frac 1 j \biggl(\sum_{i=1}^{n-1} {\binom n i} (n-i)\pder 
{u_{j-1}}{T_i} 
T_{n-i}\biggr),
\end{equation*}
implying
\begin{equation}\label{orig4.5}
\pder {u_j} T=\frac 1 j \biggl(t \pder t T\biggr) \biggl(\pder {u_{j-1}} 
T\biggr).
\end{equation}
We can then conclude from  \eqref{orig4.4} and  \eqref{orig4.5} that (cf.~\cite[(3.4)]{BB})
\begin{equation}\label{orig4.6}
\pder {u_j} T = \frac 1 {j!} \biggl(t\pder t T\biggr)^j.
\end{equation}
This equation describes a correspondence between hypertrees rooted at an unlabeled $j$-edge
and sets of $j$ hypertrees each rooted at a vertex. Husimi 
\cite{H} was the first to obtain \eqref{orig4.6}.

We now return to \eqref{orig4.3}, the definition of $T$. We apply the operator 
$u_j\pder {u_j} {}$ to both sides; note that this will count hypertrees rooted at a labeled $j$-edge. Writing $\sum_{\lambda}\tau_{n,\lambda} 
\ul$ for $T_n$, we get
\begin{equation*}
\begin{split}
u_j\pder {u_j}T&=\sum_{n=1}^{\infty} \frac {t^n} {n!} u_j \pder {u_j}{T_n}\\
&=\sum_{n=1}^{\infty} \frac {t^n} {n!} u_j \pder {u_j}{} \biggl( 
\sum_{\lambda}\tau_{n,\lambda} 
\ul\biggr)\\
&=\sum_{n=1}^{\infty} \frac {t^n} {n!} \sum_{\lambda}\tau_{n,\lambda}\, 
\lambda_j \,\ul.
\end{split}
\end{equation*}
Multiplying both sides by $j-1$ and then summing on $j$ yields
\begin{equation*}
\begin{split}
\sum_{j=2}^{\infty} (j-1)u_j \pder {u_j} T&=\sum_{j=2}^{\infty} (j-1) 
\sum_{n=1}^{\infty} \frac {t^n} {n!} \sum_{\lambda}\tau_{n,\lambda} \lambda_j 
\ul\\
&=\sum_{n=1}^{\infty} \frac {t^n} {n!} \sum_{\lambda}\tau_{n,\lambda} \ul 
\sum_{j=2}^n (j-1)\lambda_j\\
&=\sum_{n=1}^{\infty}  \frac {t^n} {n!}\sum_{\lambda}\tau_{n,\lambda} \ul (n-
1)\\
&=\sum_{n=1}^{\infty} (n-1) \frac {t^n} {n!} T_n.
\end{split}
\end{equation*}
In the above, the third equality follows from the second because $\ll j$ is 
the number of 
$j$-edges in a hypertree, and the edge magnitude of a hypertree is $n-1$. 
We conclude that (cf.~\cite[(3.5)]{BB})
\begin{equation}\label{orig4.7}
\sum_{j=2}^{\infty} (j-1)u_j \pder {u_j} T=t\pder t T-T.
\end{equation}
This equation describes two ways to count each hypertree on $n$ vertices with 
multiplicity $n-1$. It 
is clear that the right hand side does this. The terms on the left side count every 
hypertree rooted at a labeled $j$-edge $j-1$ times. But since the edge magnitude of a 
hypertree on $[n]$ is 
$n-1$, the left 
side counts every hypertree $n-1$ times. 

We now define $R$ to be
\begin{equation*}
R = t\pder t T.
\end{equation*}
In  \cite[(3.9)]{BB}, the expression $w$ corresponds to $R$, which is the 
exponential generating function for hypertrees 
rooted at a labeled vertex, counting hypertrees by weight and number of edges.
We shall refer to the objects counted by $R$ as {\it rooted hypertrees\/}.
Using this definition, and using \eqref{orig4.6} and \eqref{orig4.7}, we obtain
\begin{equation}\label{orig4.8}
T=R-\sum_{j=2}^{\infty}(j-1) u_j \frac 1 {j!} R^{j},
\end{equation}
which was also first derived by Husimi \cite{H}.

\figureform{bbfig.2}{(a) A hypertree $H$ }
\figureform{bbfig.3}{(b) Decomposition of $H$ into sets of hypertrees}
\centerline{Figure 4: Decomposition of \eqref{orig4.9} }
\bigskip\bigskip

We now obtain a functional equation for $R$, using a slightly different path from that in \cite{H}. Differentiating both sides of \eqref{orig4.8} with respect to $t$, we get
\begin{equation*}
\begin{split}
\frac {R} {t} =\pder t T
&=  \pder t {R} - \sum_{j=2}^{\infty} (j-1)u_j \frac 1 {j!} j \,R^{j-1}\, \pder t {R},
\end{split}
\end{equation*}
so 
\begin{equation}
\frac {dt} t = \frac {\partial R}{R} - 
 \sum_{j=2}^{\infty} u_j \frac 1 {(j-2)!}  R^{j-2}\, \partial R.
\end{equation}
Integrating this yields
\begin{equation}
\frac {R} t = \exp\biggl(\sum_{j=1}^{\infty} u_{j+1}\frac {R^j}{j!}\biggr),
\end{equation}
where we obtain the constant of integration by noting that
$\frac{R}{t}|_{t=0}=1$
(since the 
number of rooted hypertrees on a single vertex is 1). 
We can rewrite this as (cf.~\cite{H})
\begin{equation}\label{orig4.9}
R=t\exp\biggl(\sum_{j=1}^{\infty} u_{j+1} \frac {R^j}{j!}\biggr).
\end{equation}

Equation \eqref{orig4.9} describes a way of decomposing rooted hypertrees into a set of 
other rooted 
hypertrees. Note that in a rooted hypertree, if $v_1$ and $v_2$ are 
both adjacent to the root of the original hypertree, then $v_1$ and $v_2$ cannot 
both be in an edge which does not contain the root.
Thus, if the root of a hypertree is contained in $i$ edges containing 
$j_1+1, 
j_2+1,\ldots, j_i+1$ vertices, then when we remove the root and those edges from the 
original hypertree, we are left 
with $i$ sets of hypertrees, containing $j_1$, $j_2$, $j_3, \ldots, 
j_i$ hypertrees. 
In addition, each hypertree in each set is rooted at the vertex formerly in an edge 
with the root.

This is exactly what \eqref{orig4.9} is describing. For a given $j$, the term $u_{j+1} \frac {R^j}{j!}$
corresponds to a set of $j$ rooted hypertrees and another edge of $j+1$ vertices; this extra edge consists of the roots of the $j$ hypertrees and a new vertex (counted by the leading $t$ in \eqref{orig4.9}) which becomes the root of the new hypertree. 

Figure~4 depicts a hypertree rooted at the vertex labeled 1 and, below,  the decomposition 
resulting 
from removing the root and all edges containing it. The roots of the hypertrees 
resulting from the decomposition are denoted by larger dots. The original hypertree  that is
shown is decomposed into three sets of hypertrees, indicated in the figure.

We note here that \eqref{orig4.9}  and \eqref{orig4.8} can be obtained from  \eqref{orig2} (which is \cite[(1.6)]{BB})
and \eqref{orig3} (which is \cite[(1.7)]{BB}), respectively. In \eqref{orig2} and \eqref{orig3}, we substitute $t$ for 
$y$; $R$ for $w$; $T$ for $y\phi(y)$; and set $\phi(x)$ equal to the power series 
$\sum_{i=1}^{\infty} \frac {u_{i+1} x^i} {(i+1)!}$.

To count hypertrees only according to the total number of vertices, we can set $u_j=1$ in \eqref{orig4.8} for all $j$. Let $\Tv$ be the expression 
obtained from 
$T$ by this substitution, and let $\Rv$ be the analogous expression for $R$.
From \eqref{orig4.8}, we get a simple expression 
relating $\Tv$ to $\Rv$, where each exponential generating function now counts 
hypertrees only
by number of vertices:
\begin{equation}\label{orig4.10}
\begin{split}
\Tv&=\Rv-\Rv\sum_{j=1}^{\infty} \frac {\Rv^j}{j!}+\sum_{j=2}^{\infty}\frac 
{\Rv^j}{j!}\\
&=\Rv-\Rv(e^{\Rv}-1)+(e^{\Rv}-1-\Rv)\\
&=(e^{\Rv}-1)-\Rv(e^{\Rv}-1)\\
&=(e^{\Rv}-1)(1-\Rv).
\end{split}
\end{equation}

We can understand \eqref{orig4.10} by considering the penultimate form of the equation. 
The expression
\begin{equation}\label{orig4.11}
u_2 \pder {u_2} T + u_3 \pder {u_3} T + u_4 \pder {u_4} T + \cdots
\end{equation}
counts hypertrees rooted at an edge. Thus, each unrooted hypertree $H$  is counted
$e(H)$ times in \eqref{orig4.11}.
From \eqref{orig4.6}, we see that if we replace each $u_j$ by 1 in \eqref{orig4.11}, the resulting expression  is equal to $e^{\Rv}-\Rv-1$. But each unrooted
hypertree $H$  is counted $v(H)$ times by $\Rv$,  so that each hypertree $H$  is 
counted by $e^{\Rv}-1$ with 
multiplicity $e(H)+v(H)$. On the other hand, we can decompose a hypertree rooted 
at an edge and a vertex in the rooted edge by removing the rooted edge (but no 
vertices). We are left with a hypertree rooted at the previously rooted vertex and a 
set of hypertrees each rooted at a vertex which used to be in the 
rooted edge. These objects are exactly counted by $\Rv(e^{\Rv}-1)$, which 
therefore counts (with 
multiplicity one) each hypertree rooted at an edge and a vertex in that edge. If as before we denote 
the number of 
$i$-edges of $H$ by $\ll i$, then the number of ways to root it at an edge and a vertex in that edge 
is $\sum_i i \ll i$. But
\begin{equation*}
\sum_i i \ll i = \sum_i (i-1) \ll i + \sum_i \ll i = (v(H)-1) + e(H).
\end{equation*}
Therefore, in $\Rv(e^{\Rv}-1)$, each hypertree $H$ is counted $v(H)-1 + 
e(H)$ 
times, 
and so subtracting that expression from $e^{\Rv}-1$ produces an expression in 
which every 
hypertree is counted exactly once.  This explains \eqref{orig4.10}.

\section {Application to enumeration of hypertrees}
\label{newenumer}
By Lagrange inversion, we can find an explicit formula for rooted hypertrees by weight and number of edges.
The numbers are well-known; cf.~Husimi \cite{H}, Greene and Iba \cite{GI}, and 
Kreweras 
\cite{Ka} (in 
which hypertrees are called ``dendroids'').  

Since we can write, from \eqref{orig4.9},
\begin{equation}\label{new6.1}
R=t\prod_{j=1}^{\infty} e^{u_{j+1}\frac {R^j}{j!}},
\end{equation}
we get, using Lagrange inversion (\cite[Theorem 1.2.4, p.~17]{GJ}),
\begin{equation*}
[t^n]R=\frac 1 n [t^{n-1}]\prod_{j=1}^{\infty} \biggl(1+\bigl(nu_{j+1} \frac{t^j}{j!}\bigr) +
\frac{\bigl(n u_{j+1}\frac{t^j}{j!}\bigr)^2} {2!} + \frac{\bigl(n u_{j+1}\frac{t^j} {j!}\bigr)^3} 
{3!} +\cdots\biggr).
\end{equation*}
Letting $\lambda \vdash n-1$ denote that $\lambda$ is a partition of $n-1$ and $a_i$ 
denote the number of parts of size $i$ in $\lambda$,  we calculate
\begin{equation}\label{orig4.12}
\biggl[\frac {t^n}{n!}\biggr]R= 
\sum_{\lambda \vdash n-1}
{\binom {n-1} {\ll 1,\ll 2,\ldots}}\prod_{i} 
\frac {(nu_{i+1})^{a_i}} {a_i!}.
\end{equation}

Since there are $n$ ways to root a hypertree on $n$ vertices, $\bigl[\frac {t^n}{n!}\bigr]T=\frac 1 n \bigl[\frac {t^n}{n!}\bigr]R$, so 
for hypertrees on no more than 6 vertices,
\begin{equation*}
\begin{split}
\biggl[\frac t {1!}\biggr]T&=1,\\
\biggl[\frac {t^2 }{2!}\biggr]T&=u_2,\\
\biggl[\frac {t ^3}{3!}\biggr]T&=u_3+3u_2^2,\\
\biggl[\frac {t ^4}{4!}\biggr]T&=u_4+12u_2u_3+16u_2^3,\\
\biggl[\frac {t 
^5}{5!}\biggr]T&=u_5+20u_2u_4+15u_3^2+150u_3u_2^2+125u_2^4,\\
\biggl[\frac {t ^6}{6!}\biggr]T
&=u_6+30u_2u_5+60u_4u_3+360u_4u_2^2+540u_2u_3^2+2160u_3u_1^3+1296u_2^5
\end{split}
\end{equation*}

If we set $u_j=u$ for all $j$ in \eqref{new6.1}, we can obtain the enumerator for rooted 
hypertrees on $[n]$ by number of edges. If we let $\Rbar$ be the generating function resulting from setting $u_j=u$ in $R$, then from \eqref{new6.1}, $\Rbar=te^{u(e^{\Rbar}-1)}$. However, $e^{u(e^t-1)}=\sum_{n\ge 0} \frac {t^n}{n!}\biggl(\sum_{k=0}^{n} S(n,k)u^k\biggr)$ is the generating function for Stirling numbers of the second kind (cf.~\cite[p.~34]{StanI}). Therefore, Lagrange inversion yields
\begin{equation*}
\begin{split}
\biggl[\frac {t^n}{n!}\biggr]\Rbar&=\frac 1 n n![t^{n-1}] e^{nu(e^t-1)}\\
&=\sum_{k=1}^{n-1} (nu)^k S(n-1,k),\\
\end{split}
\end{equation*}
so the number of rooted hypertrees on $[n]$ with $k$ edges is $n^k S(n-1,k)$. In particular, the total number of rooted hypertrees on $[n]$ is $\sum_{k} n^k S(n-
1,k)$.

\end{document}